\documentclass[12pt,savetrees]{amsart}
\usepackage{amssymb, hyperref}

\newcommand{\Bgp}{{\Z^\N}}

\long\def\forget#1\forgotten{}
\newcommand{\issuenumber}{35}
\newcommand{\issuemonth}{March}
\newcommand{\issueyear}{2013}

\newcommand{\Ga}{\Gamma}
\newcommand{\Op}{\mathrm{O}}

\newcommand{\ed}{
\newpage

\section{Unsolved problems from earlier issues}

\begin{issue}Is $\binom{\Omega}{\Gamma}=\binom{\Omega}{\Tau}$?\end{issue}
\begin{issue}Is $\ufin(\Op,\Omega)=\sfin(\Gamma,\Omega)$?And if not, does $\ufin(\Op,\Gamma)$ imply
$\sfin(\Gamma,\Omega)$?\end{issue}
\stepcounter{issue}\begin{issue}Does $\sone(\Omega,\Tau)$ imply $\ufin(\Gamma,\Gamma)$?\end{issue}
\begin{issue}Is $\fp=\fp^*$? (See the definition of $\fp^*$ in that issue.)\end{issue}
\begin{issue}Does there exist (in ZFC) an uncountable set satisfying $\sfin(\cB,\cB)$?\end{issue}
\stepcounter{issue}
\begin{issue}Does $X \nin \NON(\cM)$ and $Y\nin\mathsf{D}$ imply that $X\cup Y\nin \COF(\cM)$?\end{issue}
\begin{issue}[CH]Is $\split(\Lambda,\Lambda)$ preserved under finite unions?\end{issue}
\begin{issue}Is $\cov(\cM)=\fo$? (See the definition of $\fo$ in that issue.)\end{issue}
\stepcounter{issue}
\begin{issue}Could there be a Baire metric space $M$ of weight $\aleph_1$ and a partition
$\mathcal{U}$ of $M$ into $\aleph_1$ meager sets where for each ${\mathcal U}'\subset\mathcal U$,
$\bigcup {\mathcal U}'$ has the Baire property in $M$?\end{issue}
\stepcounter{issue} 
\begin{issue}Does there exist (in ZFC) a set of reals $X$ of cardinality $\fd$ such that all
finite powers of $X$ have Menger's property $\sfin(\Op,\Op)$?\end{issue}
\begin{issue}Can a Borel non-$\sigma$-compact group be generated by a Hurewicz subspace?\end{issue}
\begin{issue}[MA]Is there $X\sbst\bbR$ of cardinality continuum, satisfying $\sone(\BO,\BG)$?\end{issue}
\begin{issue}[CH]Is there a totally imperfect $X$ satisfying $\ufin(\Op,\Gamma)$
that can be mapped continuously onto $\Cantor$?\end{issue}
\begin{issue}[CH]Is there a Hurewicz $X$ such that $X^2$ is Menger but not Hurewicz?\end{issue}
\begin{issue}Does the Pytkeev property of $C_p(X)$ imply that $X$ has Menger's property?\end{issue}
\begin{issue}Does every hereditarily Hurewicz space satisfy $\sone(\BG,\BG)$?\end{issue}
\begin{issue}[CH]Is there a Rothberger-bounded $G\le\Bgp$ such that $G^2$ is not Menger-bounded?\end{issue}
\begin{issue}Let $\cW$ be the van der Waerden ideal. Are $\cW$-ultrafilters closed under products?\end{issue}
\begin{issue}Is the $\delta$-property equivalent to the $\gamma$-property $\binom{\Omega}{\Gamma}$?\end{issue}
\stepcounter{issue}\stepcounter{issue}
\general\end{document}}

\newcommand{\x}{\times}
\newcommand{\Cantor}{{\{0,1\}^\N}}

\newcommand{\fd}{\mathfrak{d}}
\newcommand{\fp}{\mathfrak{p}}\newcommand{\fs}{\mathfrak{s}}
\newcommand{\NON}{{\mathsf   {NON}}}\newcommand{\COF}{{\mathsf   {COF}}}

\newcommand{\cN}{\mathcal{N}}
\newcommand{\cM}{\mathcal{M}}
\newcommand{\cov}{\mathsf{cov}}\newcommand{\add}{\mathsf{add}}

\newcommand{\CH}{the Continuum Hypothesis}\newcommand{\bbR}{\mathbb{R}}
\newcommand{\fo}{\mathfrak{od}}

\renewcommand{\split}{\mathsf{Split}}\newcommand{\bq}{\begin{quote}}\newcommand{\eq}{\end{quote}}
\newcommand{\cB}{\mathcal{B}}\newcommand{\BG}{\cB_\Gamma}
\newcommand{\BO}{\cB_\Omega}

\newcommand{\sone}{\mathsf{S}_1}\newcommand{\sfin}{\mathsf{S}_\mathrm{fin}}
\newcommand{\ufin}{\mathsf{U}_\mathrm{fin}} 
\newcommand{\nin}{\not\in}
\newcommand{\cW}{\mathcal{W}}

\newcommand{\N}{\mathbb{N}}\newcommand{\Z}{\mathbb{Z}}
\newcommand{\sbst}{\subseteq}
\newcommand{\by}[2]{\par\hfill\emph{#1}, #2}\newcommand{\nby}[1]{\par\hfill\emph{#1}}\newcommand{\Tau}{\mathrm{T}}
\newcommand{\CE}{\textsc{CE}}
\newtheorem{thm}{Theorem}[section]\newcommand{\bthm}{\begin{thm}} \newcommand{\ethm}{\end{thm}}
\newtheorem{prop}[thm]{Proposition}\newcommand{\bprp}{\begin{prop}} \newcommand{\eprp}{\end{prop}}
\newtheorem{fact}[thm]{Fact}\newcommand{\bfct}{\begin{fact}} \newcommand{\efct}{\end{fact}}
\newtheorem{prob}[thm]{Problem}\newcommand{\bprb}{\begin{prob}} \newcommand{\eprb}{\end{prob}}
\newtheorem{lem}[thm]{Lemma}\newcommand{\blem}{\begin{lem}} \newcommand{\elem}{\end{lem}}
\newtheorem{claim}[thm]{Claim}\newcommand{\bclm}{\begin{claim}} \newcommand{\eclm}{\end{claim}}
\newtheorem{cor}[thm]{Corollary}\newcommand{\bcor}{\begin{cor}} \newcommand{\ecor}{\end{cor}}
\newtheorem{conj}[thm]{Conjecture}\newcommand{\bcnj}{\begin{conj}} \newcommand{\ecnj}{\end{conj}}
\theoremstyle{definition}\newtheorem{defn}[thm]{Definition}\newcommand{\bdfn}{\begin{defn}} \newcommand{\edfn}{\end{defn}}
\theoremstyle{remark}\newtheorem{rem}[thm]{Remark}\newcommand{\brem}{\begin{rem}} \newcommand{\erem}{\end{rem}}
\newtheorem{cnv}[thm]{Convention}\newcommand{\bcnv}{\begin{cnv}} \newcommand{\ecnv}{\end{cnv}}
\newtheorem{exam}[thm]{Example}\newcommand{\bexm}{\begin{exam}} \newcommand{\eexm}{\end{exam}}
\newtheorem{issue}{Issue}\newcommand{\bpf}{\begin{proof}} \newcommand{\epf}{\end{proof}}
\newcommand{\be}{\begin{enumerate}}\newcommand{\ee}{\end{enumerate}}\newcommand{\bi}{\begin{itemize}}
\newcommand{\ei}{\end{itemize}}
\newcommand{\general}{\small\vfill\par\noindent\hrulefill\par
\noindent\textbf{Previous issues.} The previous issues of this
bulletin are available online at\\
\url{http://front.math.ucdavis.edu/search?\&t=\%22SPM+Bulletin\%22}
\\[0.1cm]
\textbf{Contributions.} Announcements, discussions, and open problems should be emailed
to \texttt{tsaban@math.biu.ac.il}\\[0.1cm]
\textbf{Subscription.}
To receive this bulletin (free) to your e-mailbox, e-mail us.
}

\newcommand{\arXivl}[4]{\subsection{#2}{#4}\par\hfill{\arx{#1}}\par\hfill\emph{#3}}
\newcommand{\arXiv}[3]{\subsection{#2}\mbox{}\par\hfill{\arx{#1}}\par\hfill\emph{#3}}

\newcommand{\AMS}[3]{\subsection{#1}\mbox{}\par\hfill{\texttt{#3}}\par\hfill\emph{#2}}

\newcommand{\arx}[1]{\url{http://arxiv.org/abs/#1}}

\title[$\mathcal{SPM}$ Bulletin \textbf{\issuenumber} (\issuemonth{} \issueyear)]{%
$\mathcal{SPM}$ Bulletin\\[0.5cm]
Issue number \issuenumber: \issuemonth{} \issueyear{} \CE{}}

\begin{document}
\maketitle


\section{Editor's note}

This issue contains a relatively large number of abstracts of papers dealing directly with selection principles.
This is, in part, due to the stimulus created by the last SPM conference, followed by a special issue 
of Topology and its Applications that is under preparation.

Pay special attention to Section \ref{Mal} below, that announces a solution of the classic Malyhin's Problem.
This problem is closely related to the question of existence of $\gamma$-sets (one of the central themes of selection principles).
In addition to its main result, this paper establishes that the existence of Malyhin groups need not imply that of $\gamma$-sets.
Indeed, the result is quite flexible and may have additional consequences in the realm of SPM. 

\medskip

With best regards,

\by{Boaz Tsaban}{tsaban@math.biu.ac.il}

\hfill \texttt{http://www.cs.biu.ac.il/\~{}tsaban}

\section{Long announcements}

\arXivl{1210.2118}
{A characterization of the Menger property by means of ultrafilter
  convergence}
{Paolo Lipparini}
{We characterize various Menger-related properties by means of ultrafilter
convergence, and discuss their behavior with respect to products.}

\arXiv{1210.2120}
{Topological spaces compact with respect to a set of filters}
{Paolo Lipparini}
{If $\mathcal P$ is a family of filters over some set $I$, a topological space
$X$ is sequencewise $\mathcal P$-compact if, for every $I$-indexed sequence of
elements of $X$, there is $F \in \mathcal P$ such that the sequence has an
$F$-limit point. Countable compactness, sequential compactness, initial
$\kappa$-compactness, $[ \lambda ,\mu]$-compactness, the Menger and Rothberger
properties can all be expressed in terms of sequencewise $\mathcal
P$-compactness, for appropriate choices of $\mathcal P$.

  We show that sequencewise $\mathcal P$-compactness is preserved under taking
products if and only if there is a filter $F \in \mathcal P$ such that
sequencewise $\mathcal P$-compactness is equivalent to $F$-compactness. If this
is the case, and there exists a sequencewise $\mathcal P$-compact $T_1$
topological space with more than one point, then $F$ is necessarily an
ultrafilter.}

\arXivl{1210.4986}
{Comparing weak versions of separability}
{Daniel T. Soukup, Lajos Soukup, Santi Spadaro}
{Our aim is to investigate spaces with sigma-discrete and meager dense sets,
as well as selective versions of these properties. We construct numerous
examples to point out the differences between these classes while answering
questions of Tkachuk [30], Hutchinson [17] and the authors of [8].}

\arXivl{1210.8010}
{Productively Lindel\"of and indestructibly Lindel\"of spaces}
{Haosui Duanmu, Franklin D. Tall, Lyubomyr Zdomskyy}
{There has recently been considerable interest in productively Lindel\"of
spaces, i.e. spaces such that their product with every Lindel\"of space is
Lindel\"of. Here we make several related remarks about such spaces.
Indestructible Lindel\"of spaces, i.e.\ spaces that remain Lindel\"of in every
countably closed forcing extension, were introduced by Tall in 1995. Their
connection with topological games and selection principles was explored by
Scheepers and Tall in 2010. We find further connections here.}

\arXivl{1211.1719}
{Indestructibility of compact spaces}
{Rodrigo R. Dias and Franklin D. Tall}
{In this article we investigate which compact spaces remain compact under
countably closed forcing. We prove that, assuming the Continuum Hypothesis, the
natural generalizations to $\omega_1$-sequences of the selection principle and
topological game versions of the Rothberger property are not equivalent, even
for compact spaces. We also show that Tall and Usuba's "$\aleph_1$-Borel
Conjecture" is equiconsistent with the existence of an inaccessible cardinal.}

\arXivl{1211.3581}
{Some observations on compact indestructible spaces}
{Angelo Bella}
{Inspired by a recent work of Dias and Tall, we show that a compact
indestructible space is sequentially compact. We also prove that a Lindelof
Hausdorff indestructible space has the finite derived set property and a
compact Hausdorff indestructible space is pseudoradial.}

\arXivl{1211.2764}
{Reflecting Lindel\"of and converging $\omega_1$-sequences}
{Alan Dow and Klaas Pieter Hart}
{We deal with a conjectured dichotomy for compact Hausdorff spaces: each such
space contains a non-trivial converging omega-sequence or a non-trivial
converging $\omega_1$-sequence. We establish that this dichotomy holds in a
variety of models; these include the Cohen models, the random real models and
any model obtained from a model of CH by an iteration of property K posets. In
fact in these models every compact Hausdorff space without non-trivial
converging $\omega_1$-sequences is first-countable and, in addition, has many
$\aleph_1$-sized Lindel\"of subspaces. As a corollary we find that in these models
all compact Hausdorff spaces with a small diagonal are metrizable.}

\arXivl{1212.0589}
{Selections, games and metrisability of manifolds}
{David Gauld}
{In this note we relate some selection principles to metrisability and
separability of a manifold. In particular we show that $\sf{S}_{fin}(\mathcal
K,\mathcal O)$, $\sf{S}_{fin}(\Omega,\Omega)$ and
$\sf{S}_{fin}(\Lambda,\Lambda)$ are each equivalent to metrisability for a
manifold, while $\sf{S}_1(\mathcal D,\mathcal D)$ is equivalent to separability
for a manifold.}

\arXivl{1212.5724}
{Infinite games and cardinal properties of topological spaces}
{Angelo Bella and Santi Spadaro}
{Inspired by work of Scheepers and Tall, we use properties defined by
topological games to provide bounds for the cardinality of topological spaces.
We obtain a partial answer to an old question of Bell, Ginsburg and Woods
regarding the cardinality of weakly Lindel\"of first-countable regular spaces
and answer a question recently asked by Babinkostova, Pansera and Scheepers. In
the second part of the paper we study a game-theoretic version of cellularity,
a special case of which has been introduced by Aurichi. We obtain a
game-theoretic proof of Shapirovskii's bound for the number of regular open
sets in an (almost) regular space and give a partial answer to a natural
question about the productivity of a game strengthening of the countable chain
condition that was introduced by Aurichi. As a final application of our results
we prove that the Hajnal-Juh\'asz bound for the cardinality of a
first-countable ccc Hausdorff space is true for almost regular (non-Hausdorff)
spaces.}

\arXivl{1212.6122}
{Weak covering properties and selection principles}
{L. Babinkostova, B. A. Pansera and M. Scheepers}
{No convenient internal characterization of spaces that are productively
Lindelof is known. Perhaps the best general result known is Alster's internal
characterization, under the Continuum Hypothesis, of productively Lindelof
spaces which have a basis of cardinality at most $\aleph_1$. It turns out that
topological spaces having Alster's property are also productively weakly
Lindelof. The weakly Lindelof spaces form a much larger class of spaces than
the Lindelof spaces. In many instances spaces having Alster's property satisfy
a seemingly stronger version of Alster's property and consequently are
productively X, where X is a covering property stronger than the Lindelof
property. This paper examines the question: When is it the case that a space
that is productively X is also productively Y, where X and Y are covering
properties related to the Lindelof property.}

\arXiv{1301.3484}
{Asymptotic dimension, decomposition complexity, and Haver's property C}
{Alexander Dranishnikov and Michael Zarichnyi}
{The notion of the decomposition complexity was introduced in [GTY] using a
game theoretical approach. We introduce a notion of straight decomposition
complexity and compare it with the original as well with the asymptotic
property C. Then we define a game theoretical analog of Haver's property C in
the classical dimension theory and compare it with the original.}

\subsection{Malykhin's Problem}\label{Mal}
We construct a model of ZFC where every separable Fr\'echet--Urysohn 
group is metrizable. This solves a 1978 problem of V. I. Malykhin.

\url{www.matmor.unam.mx/~michael/preprints_files/Frechet-malykhin.pdf}

\nby{Michael Hru\v{s}\'ak and Ulises Ariet Ramos-Garc\'ia}

\arXiv{1302.5287}
{A new class of spaces with all finite powers Lindelof}
{Natasha May, Santi Spadaro and Paul Szeptycki}
{We consider a new class of open covers and classes of spaces defined from
them, called "iota spaces". We explore their relationship with epsilon-spaces
(that is, spaces having all finite powers Lindelof) and countable network
weight. An example of a hereditarily epsilon-space whose square is not
hereditarily Lindelof is provided.}

\arXiv{1302.5658}
{Topologically invariant $\sigma$-ideals on the Hilbert cube}
{Taras Banakh, Michal Morayne, Robert Ralowski, Szymon Zeberski}
{We study and classify topologically invariant $\sigma$-ideals with a Borel
base on the Hilbert cube and evaluate their cardinal characteristics. One of
the results of this paper solves (positively) a known problem whether the
minimal cardinalities of the families of Cantor sets covering the unit interval
and the Hilbert cube are the same.}

\arXiv{1303.0815}
{Topological spaces compact with respect to a set of filters. II}
{Paolo Lipparini}
{If $\mathcal P$ is a family of filters over some set $I$, a topological space
$X$ is sequencewise $\mathcal P$ compact if, for every $I$-indexed sequence of
elements of $X$, there is $F \in \mathcal P$ such that the sequence has an
$F$-limit point. As recalled in Part I, countable compactness, sequential
compactness, initial $\kappa$-compactness, $[ \lambda ,\mu]$-compactness, the
Menger and Rothberger properties can all be expressed in terms of sequencewise
$\mathcal P$ compactness, for appropriate choices of $\mathcal P$.
  We show that a product of topological spaces is sequencewise $\mathcal P$
compact if and only if so is any subproduct with $\leq |\mathcal P|$ factors.
In the special case of sequential compactness, we get a better bound: a product
is sequentially compact if and only if all subproducts by $\leq \fs$ factors
are sequentially compact.}

\arXiv{1303.3597}
{Selective covering properties of product spaces}
{Arnold W. Miller, Boaz Tsaban, Lyubomyr Zdomskyy}
{We study the preservation of selective covering properties, including classic
ones introduced by Menger, Hurewicz, Rothberger, Gerlits and Nagy, and others,
under products with some major families of concentrated sets of reals.

  Our methods include the projection method introduced by the authors in an
earlier work, as well as several new methods. Some special consequences of our
main results are (definitions provided in the paper): 
\be
  \item Every product of a concentrated space with a Hurewicz $\sone(\Ga,\Op)$
space satisfies $\sone(\Ga,\Op)$. On the other hand, assuming \CH{}, for each
Sierpi\'nski set $S$ there is a Luzin set $L$ such that $L\x S$ can be mapped
onto the real line by a Borel function.
  \item Assuming Semifilter Trichotomy, every concentrated space is
productively Menger and productively Rothberger.
  \item Every scale set is productively Hurewicz, productively Menger,
productively Scheepers, and productively Gerlits--Nagy.
  \item Assuming $\fd=\aleph_1$, every productively Lindel\"of space is
productively Hurewicz, productively Menger, and productively Scheepers. 
\ee
  A notorious open problem asks whether the additivity of Rothberger's property
may be strictly greater than $\add(\cN)$, the additivity of the ideal of
Lebesgue-null sets of reals. We obtain a positive answer, modulo the
consistency of Semifilter Trichotomy with $\add(\cN)<\cov(\cM)$.

  Our results improve upon and unify a number of results, established earlier
by many authors.}

\section{Short announcements}\label{RA}

\arXiv{1210.2121}
{Productivity of $[\mu, \lambda ]$-compactness}
{Paolo Lipparini}

\AMS{Almost Souslin Kurepa trees}
{Mohammad Golshani}
{http://www.ams.org/journal-getitem?pii=S0002-9939-2012-11461-3}

\arXiv{1212.5725}
{On two topological cardinal invariants of an order-theoretic flavour}
{Santi Spadaro}

\arXiv{1212.5726}
{P-spaces and the Volterra property}
{Santi Spadaro}

\arXiv{1301.5297}
{Hereditarily supercompact spaces}
{Taras Banakh, Zdzislaw Kosztolowicz, Slawomir Turek}

\arXiv{1303.3600}
{Hindman's Coloring Theorem in arbitrary semigroups}
{Gili Golan, Boaz Tsaban}

\ed